\documentclass[12pt,reqno]{article}      
\def\hybrid{\topmargin      0pt
\oddsidemargin 0pt 
\headheight 0pt \headsep 0pt 
\textwidth 165true mm  
\textheight 231true mm 
\marginparwidth 0.0in \parskip  0pt plus 1pt \jot = 1.5ex}

\usepackage{amssymb}  
\usepackage{amsmath}  
\usepackage{amsthm}  
\hybrid
\usepackage{amssymb}     
\usepackage{amsmath}       
\usepackage{amsthm}
\newcommand{\be}[1]{\begin{eqnarray#1}}
\newcommand{\ee}[1]{\end{eqnarray#1}} 
\newtheorem{thm}{Theorem}[section]
\newtheorem{propn}[thm]{Proposition}      
\newtheorem{lemma}[thm]{Lemma}
\newtheorem{corollary}[thm]{Corollary}

\theoremstyle{definition}
\newtheorem{remark}[thm]{Remark}
\newtheorem{definition}[thm]{Definition}
\newtheorem{example}[thm]{Example}

\newcommand{\tp}{\otimes}

\newcommand{\braid}[2]{{#1}$\lower4pt\hbox{$\tp\atop\raise4pt
            \hbox{$\scriptscriptstyle\Ru $}$}${#2}}
\newcommand{\twist}[2]{{#1}${\,\scriptscriptstyle \Ru}\atop\raise9pt\hbox{$\scriptstyle\tp$} ${#2}}
\newcommand{\twistF}[2]{{#1}${\,\scriptscriptstyle \F}\atop\raise9pt\hbox{$\scriptstyle\tp$} ${#2}}

\newcommand{\A}{\mathcal{A}}
\newcommand{\B}{\mathcal{B}}                
\newcommand{\M}{\mathcal{M}}
\newcommand{\Ha}{\mathcal{H}}              
\newcommand{\Ru}{\mathcal{R}}
\newcommand{\Q}{\mathcal{Q}}
\newcommand{\U}{\mathcal{U}}                
\newcommand{\F}{\mathcal{F}}
\newcommand{\C}{\mathbb{C}}                  
                  
\newcommand{\K}{\mathbb{K}}

\newcommand{\ad}{\mathrm{ad}}              
\newcommand{\La}{\mathcal{L}}
\newcommand{\Ta}{\mathcal{T}}

\newcommand{\g}{\mathfrak{g}}

\renewcommand{\[}{{[\![}}

\newcommand{\D}{\mbox{D}}

\newcommand{\ve}{\varepsilon}                   

\newcommand{\n}{\nonumber          }           
\newcommand{\al}{\alpha}                         
\newcommand{\bt}{\beta}
                     
\newcommand{\id}{\mbox{id}}
\newcommand{\End}{\mathrm{End}}

\newcommand{\T}{\mathrm{T}}

\newcommand{\gm}{\gamma}                           

\newcommand{\tr}{\triangleright}         
\newcommand{\tl}{\triangleleft}

\begin{document}
\title{Reflection Equation, Twist, and Equivariant
Quantization}
\author{J. Donin\footnote{This research is partially supported 
by the Israel Academy of Sciences grant no. 8007/99-01.}
\hspace{3pt} and A. Mudrov\\[0.1in]
{\small
Department of Mathematics, Bar-Ilan University, 52900 Ramat-Gan, Israel}
}
\date{}
\maketitle
\begin{abstract}
We prove that the reflection equation (RE) algebra $\La_R$ associated
with a finite dimensional representation of a quasitriangular
Hopf algebra $\Ha$ is twist-equivalent
to the corresponding Faddeev-Reshetikhin-Takhtajan (FRT) algebra.
We show that $\La_R$ is a module algebra over the twisted
tensor square \twist{$\Ha$}{$\Ha$} and the double $\D(\Ha)$.
We define FRT- and RE-type algebras and apply them to 
 the problem of equivariant quantization
on Lie groups and matrix spaces.
\end{abstract}

\section{Introduction.}
Let $\Ha$ be a quasitriangular Hopf algebra with the universal R-matrix $\Ru$ .
Let $V$ be the  space of its finite dimensional representation and
$R$ the image of $\Ru$ in $\End^{\tp2}(V)$. We study relations between
two algebras naturally arising in this context, the 
Faddeev-Reshetikhin-Takhtajan (FRT) algebra $\Ta_R$ and the so-called
reflection equation (RE) algebra  $\La_R$. They are both  
quotients of the tensor algebra $\T\bigl(\End^*(V)\bigr)$ by 
quadratic relations and admit certain Hopf algebra actions.
So  $\Ta_R$ is endowed with the structure of a bimodule
algebra over $\Ha$ that is naturally extended from the bimodule
structure on $\End^*(V)$.  The algebra $\La_R$ is a left $\Ha$-module
algebra with the action extended from the coadjoint representation 
on $\End^*(V)$. 

Our first result is that the RE algebra  $\La_R$ is a module algebra
not only over $\Ha$ but over the twisted tensor square \twist{$\Ha$}{$\Ha$}.
The latter is the twist of the ordinary tensor product $\Ha\tp \Ha$ with
the universal R-matrix as a twisting cocycle, 
$\Ru_{23}\in (\Ha\tp \Ha)\tp (\Ha\tp \Ha)$. The action of $\Ha$ on
$\La_R$ is induced by the Hopf algebra embedding $\Ha\to \mbox{\twist{$\Ha$}{$\Ha$}}$
via the coproduct.

As a corollary, we obtain that  $\La_R$ is a module algebra over
the coopposite dual $\Ha^{*op}$ because there exists a Hopf homomorphism
$\Ha^{*op}\to \mbox{\twist{$\Ha$}{$\Ha$}}$. Since the Hopf algebra
homomorphisms from $\Ha$ and $\Ha^{*op}$ to \twist{$\Ha$}{$\Ha$} can be
extended to a Hopf algebra homomorphism from the double $\D(\Ha)$ to 
$\mbox{\twist{$\Ha$}{$\Ha$}}$, we obtain that $\La_R$ 
is a module algebra over $\D(\Ha)$.

Our second result is that $\La_R$ is a twist of $\Ta_R$ as 
a module algebra. The $\Ha$-bimodule $\Ta_R$ can be considered
as a left $\Ha^{op}\tp\Ha$-module. The twist from $\Ha^{op}\tp\Ha$
to \twist{$\Ha$}{$\Ha$} is performed via the cocycle $\Ru_{13}\Ru_{23}$,
where the first transformation via $\Ru_{13}$ converts the first tensor factor
$\Ha^{op}$ to $\Ha$ while the second twist  via $\Ru_{23}$ makes the ordinary tensor
square $\Ha^{\tp2}$  the twisted one.

The algebra $\Ta_R$ is commutative in the category of $\Ha$-bimodules.
Using this fact, we prove that $\La_R$ is commutative in the category of 
\twist{$\Ha$}{$\Ha$}-modules. In general, we prove twist-equivalence
between the classes of FRT- and RE-type algebras, which we define to be
commutative algebras in the categories of $\Ha$-bimodules and
\twist{$\Ha$}{$\Ha$}-modules, respectively.

In particular, we introduce the RE dual algebra $\tilde \Ha^*$ as an RE-type algebra 
that is twist-equivalent to the FRT-type algebra $\Ha^*$ and we show that 
it coincides with the braided Hopf algebra of Majid.

We study coactions  on  $\tilde \Ha^*$ of the Hopf algebras 
$(\mbox{\twist{$\Ha$}{$\Ha$}})^*$, $\Ha^*$, and the opposite Hopf algebra 
$\Ha_{op}$. We deduce 
properties of the $\Ha$-equivariant homomorphism $\tilde \Ha^*\to \Ha$, 
$\xi\to\langle\xi,\Q_1 \rangle \Q_2$, 
where $\Q=\Ru_{21}\Ru$, using the coalgebra structure over $\Ha_{op}$.

We apply our construction to the deformation quantization on 
Lie groups and matrix spaces. In particular, we show that
the algebra $\tilde \U^*_h(\g)$, where $\g$ is a semisimple 
Lie algebra, is the $\U_h(\g)$-equivariant quantization of a special Poisson structure 
on the corresponding Lie group, the RE bracket. 
It is known that the quotient 
of $\Ta_R$ by torsion
is the quantization on the cone $\End^\Omega(V)$ of matrices
whose tensor square commutes with the image $\Omega$ of the 
split Casimir, the invariant symmetric element from $\g^{\tp2}$.
As an implication of the twist-equivalence between $\La_R$ and $\Ta_R$, 
we find that the  quotient of $\La_R$ by torsion is the 
\twist{$\U_h(\g)$}{$\U_h(\g)$}-equivariant 
quantization on $\End^\Omega(V)$.

The setup of the paper is as follows. Section \ref{sQHA} contains basic 
facts about quasitriangular Hopf algebras, the twist transformation, and
the relation between the double and twisted tensor square. 
Section \ref{sMHA} recalls what are modules and comodules
over Hopf algebras. 
The relations between the FRT and RE algebra and their implications
are studied in Section~ \ref{FRT-RE-type}. 
 Section \ref{sAQGS} is devoted to applications to the equivariant
 deformation quantization on Lie groups and matrix spaces.
\newpage
\section{Quasitriangular Hopf algebras.}
\label{sQHA}
\subsection{Definitions and elementary properties.}
\label{ssDEP}
In this subsection, we recall basic definitions of the quasitriangular 
Hopf algebra theory, \cite{Dr1}. For a detailed exposition,
the reader can consult to \cite{Mj}.
Let $\K$ be a commutative algebra over a field
of zero characteristic. Let $\Ha$ be a Hopf algebra over $\K$,
with the coproduct $\Delta\colon \Ha \to \Ha\tp \Ha$, counit 
$\ve \colon \Ha \to \K$, and  antipode $\gm  \colon\Ha \to \Ha$.
Throughout the paper, we adopt the  standard notation with implicit summation
in order to explicate factors of tensor objects, e.g., we write
$\Phi = \Phi_1\tp \ldots \tp \Phi_k$ for an element $\Phi \in \Ha^{\tp k}$. 
For the coproduct, we use the symbolic Sweedler notation, 
$\Delta(x)=x_{(1)}\tp x_{(2)}$, $x\in \Ha$. 

A Hopf algebra $\Ha$ is called {\em quasitriangular}
if there is an element $\Ru \in\Ha^{\tp2}$, the universal R-matrix,
such that 
\be{}
(\Delta \tp id)(\Ru)=\Ru_{13}\Ru_{23}, \quad (id \tp\Delta)(\Ru)=\Ru_{13}\Ru_{12}
\label{bch}
\ee{}
and, for any $x\in \Ha$,
\be{}
\quad\Ru\Delta(x)=\Delta^{op}(x)\Ru.
\label{flip}
\ee{}
The subscripts in (\ref{bch}) specify the way of embedding $\Ha^{\tp 2}$ into $\Ha^{\tp 3}$, namely, 
$\Ru_{12}=\Ru_1\tp \Ru_2\tp 1$, $\Ru_{13}=\Ru_1\tp 1\tp \Ru_2$, and
$\Ru_{23}=1\tp \Ru_1 \tp \Ru_2$. The opposite coproduct $\Delta^{op}$ is the 
composition of $\Delta$ and the flip operator  $\tau$ on $\Ha^{\tp 2}$.

The following identities are implications
of defining relations  (\ref{bch}) and (\ref{flip}):
\be{}
(\ve\tp \id)(\Ru) = 1\tp 1,& & (\id\tp\ve)(\Ru)=1\tp 1,\\
(\gamma\tp \id)(\Ru) = \Ru^{-1}, &\quad&  (\id \tp \gamma)(\Ru^{-1}) = \Ru.
\label{urm-antipode}
\ee{}
Also, the Yang-Baxter equation in $\Ha^{\tp 3}$, 
\be{}
\Ru_{12}\Ru_{13}\Ru_{23}&=&\Ru_{23}\Ru_{13}\Ru_{12},
\label{YBE}
\ee{}
follows from (\ref{bch}) and (\ref{flip}).

There is an alternative quasitriangular 
structure\footnote{They may coincide.} 
on $\Ha$
with the universal R-matrix $\Ru^{-1}_{21}=\tau(R^{-1})$. Obviously,
it fulfills conditions (\ref{bch}) and (\ref{flip}).

The opposite Hopf algebra $\Ha_{op}$ is endowed with the opposite
multiplication $m_{op}=m\circ\tau$, where $m$ is the original one.
Similarly, the coopposite Hopf algebra $\Ha^{op}$ has
the coproduct $\Delta^{op}=\tau\circ\Delta$. The antipode of $\Ha$ 
is an isomorphism between $\Ha$ and $\Ha^{op}_{op}$,
being an anti-algebra and anti-coalgebra  map. Important
for our exposition is that we may also treat it as an isomorphism
between $\Ha_{op}$ and $\Ha^{op}$.

The dual Hopf algebra $\Ha^*$ is spanned by matrix coefficients of all 
finite dimensional representations\footnote{We assume that the supply 
of representations is large enough to separate elements of $\Ha$.}  
of $\Ha$. 
There are two remarkable maps from $\Ha^{*op}$ to $\Ha$ defined
via the universal R-matrix:
\be{}
\Ru^\pm(\eta) = \langle\eta,\Ru^\pm_{1}\rangle \Ru^\pm_{2}, 
\quad \eta\in \Ha^{*op},\quad \mbox{where}\quad
\Ru^+ = \Ru\quad \mbox{and }\quad\Ru^- = \Ru^{-1}_{21}.
\label{R+-}
\ee{}
It follows from (\ref{bch}) that they are Hopf algebra homomorphisms
$\Ha^{*op}\to \Ha$
(as was already mentioned, the element $\Ru^{-1}_{21}$ satisfies (\ref{bch}) as well).

The tensor product $\A\tp\B$  of two Hopf algebras $\A$ and $\B$
is a Hopf algebra with the multiplication 
\be{}
\label{tpm}
(a_1\tp b_1) (a_2\tp b_2)=a_1 a_2\tp b_1 b_2,
\quad a_i\in \A, \>b_i\in \B, \>i=1,2,
\ee{}
coproduct 
\be{}
\label{tpc}
\Delta(a\tp b)=\bigl(a_{(1)}\tp b_{(1)}\bigr)\tp \bigl(a_{(2)}\tp b_{(2)}\bigr),
\quad a\in \A, \>b\in \B,
\ee{}
counit $\ve=\ve_\A\tp \ve_\B$, and antipode $\gm=\gm_\A\tp \gm_\B$.
If $\A$ and $\B$ are both quasitriangular, so is $\A\tp \B$. Its 
universal R-matrix is
\be{}
\Ru_{\A\tp\B}=\Ru_\A\Ru_\B,
\label{R_AB}
\ee{}
where $\Ru_\A$ and $\Ru_\B$ are R-matrices of $\A$ and $\B$ naturally embedded
in $(\A\tp\B)^{\tp2}$.
\begin{remark}
In the infinite dimensional case, we assume that all algebras are complete
in some topology. It may be, for example,  the $h$-adic 
topology in the case $\K=\C[[h]]$. All tensor products are assumed
to be completed.
\end{remark}
\subsection{Twist of Hopf algebras.}
\label{ssTHA}
In this subsection, we collect several facts and examples, which
will be essential for our further exposition, concerning the 
twist transformation of Hopf algebras, \cite{Dr2}. 
Let $\F$ be an invertible element from $\Ha\otimes\Ha$ satisfying the 
cocycle constraint
\be{}
(\Delta \tp \id)(\F)\F_{12} = (\id\tp \Delta)(\F)\F_{23},
\label{cc}
\ee{}
with the normalizing condition
$(\ve \tp \id)(\F) =1\tp 1 = (\id\tp \ve)(\F)$.
There exists a new Hopf algebra structure $\tilde \Ha$ on $\Ha$ with the same 
multiplication and counit but the  "twisted" coproduct
\be{}
\label{tcopr}
\tilde \Delta(x)=\F^{-1}\Delta(x)\F, \quad x\in \Ha,
\ee{}
and antipode 
$$\tilde \gm(x) = u^{-1}\gm(h)u, \quad \mbox{where}\quad
u=\gm(\F_1)\F_2 \in \Ha.$$
Condition (\ref{cc}) ensures the coproduct $\tilde \Delta$ being
coassociative.
The Hopf algebra $\tilde \Ha$ is quasitriangular, 
provided so is $\Ha$. Its universal R-matrix
is expressed through the old one and the twisting cocycle:
\be{}
\label{Rtwisted}
\tilde \Ru = \F^{-1}_{21}\Ru \F.
\ee{}
We call algebras $\Ha$ and $\tilde \Ha$ {\em twist-equivalent}
and use the notation $\tilde\Ha\stackrel{\F}{\sim}\Ha$ or
simply $\tilde\Ha\sim\Ha$ when the exact form of $\F$ is
clear from the context. 
Obviously, 
$\Ha_2\stackrel{\F_1}{\sim}\Ha_1$ and
$\Ha_3\stackrel{\F_2}{\sim}\Ha_2$ imply
$\Ha_3\stackrel{\F_1\F_2}{\sim}\Ha_1$. Also,
if $\Ha_2\stackrel{\F}{\sim}\Ha_1$, then
$\Ha_1\stackrel{\F^{-1}}{\sim}\Ha_2$.
Essential for us will be the following examples.
\begin{example}[Coopposite Hopf algebra]
\label{COP}
Given a quasitriangular Hopf algebra $\Ha$ its coopposite algebra
$\Ha^{op}$ can be  obtained by twist with $\F=\Ru^{-1}$, cf. equation (\ref{flip}). 
The cocycle condition (\ref{cc}) follows from the Yang-Baxter equation 
(\ref{YBE}). 
\end{example}
\begin{example}[Twisted tensor product]
\label{TTP}
Let $\Ha$ be a tensor product $\Ha=\A\tp \B$ of two Hopf algebras
with multiplication (\ref{tpm}) and coproduct (\ref{tpc}).
An element $\F\in \B\tp \A$ may be viewed as that from
$\Ha\tp\Ha$  via the embedding
$(1\tp\B)\tp (\A\tp 1)\subset(\A\tp \B)\tp (\A\tp \B)$.
If $\F$ satisfies the identities
\be{}
(\Delta_{\B}\tp\id)(\F)=\F_{13}\F_{23}\in \B\tp\B\tp\A
,\quad
(\id\tp\Delta_\A)(\F)=\F_{13}\F_{12}\in \B\tp\A\tp\A,
\label{ttp}
\ee{}
it also fulfills the cocycle condition (\ref{cc}), 
\cite{RS}.
\begin{definition}
{\em Twisted tensor product} \twistF{$\A$}{$\B$} of two Hopf 
algebras is the twist of $\A\tp \B$ with a cocycle $\F$ satisfying 
(\ref{ttp}).
\end{definition}
\noindent
An immediate corollary of conditions (\ref{ttp}) is that the
evaluation maps 
$\id \tp \ve_\B\colon \mbox{\twistF{$\A$}{$\B$}}\to \A$
and 
$\ve_\A\tp \id \colon \mbox{\twistF{$\A$}{$\B$}}\to \B$
are Hopf ones. Note that, contrary to the ordinary tensor product,
the embeddings of $\A$ and $\B$ into \twistF{$\A$}{$\B$}
are algebra but not coalgebra maps.

An important example of the twisted tensor product is when
$\A=\B=\Ha$ is a quasitriangular Hopf algebra and $\F=\Ru$. 
Condition  (\ref{ttp}) then holds because of  (\ref{bch}).
This particular case is called {\em twisted tensor square} 
of a quasitriangular Hopf algebra and denoted
\twist{$\Ha$}{$\Ha$}.
It is convenient for our exposition to take
the universal R-matrix $\Ru^-_{13}\Ru^+_{24}$
for the ordinary tensor product $\Ha\tp\Ha$, see (\ref{R+-}) and (\ref{R_AB}). 
Then formula (\ref{Rtwisted}) 
gives the universal R-matrix of \twist{$\Ha$}{$\Ha$}:
\be{}
\label{ttpRm}
\Ru^{-1}_{41}\Ru^-_{13}\Ru^+_{24} \Ru_{23}=
\Ru^{-1}_{41}\Ru^{-1}_{31}\Ru _{24} \Ru_{23}
\in 
(\mbox{\twist{$\Ha$}{$\Ha$}})\tp(\mbox{\twist{$\Ha$}{$\Ha$}}).
\ee{}
\end{example}
\subsection{Drinfeld's double and twisted tensor square.}
\label{ssDTTP}
This subsection is devoted to a relation between the double $\D(\Ha)$ 
of a quasitriangular Hopf algebra $\Ha$ and its twisted tensor square
\twist{$\Ha$}{$\Ha$}. 
A particular case of the so called 
factorizable Hopf algebras was considered in \cite{RS}. For those algebras, 
the double $\D(\Ha)$ is isomorphic to  $\mbox{\twist{$\Ha$}{$\Ha$}}$. In general, 
there is a Hopf algebra homomorphism from $\D(\Ha)$ to  $\mbox{\twist{$\Ha$}{$\Ha$}}$.

As a coalgebra, the double $\D(\Ha)$ coincides with  the tensor product
$\Ha\tp \Ha^{*op}$. Both of $\Ha$ and $\Ha^{*op}$ are embedded in $\D(\Ha)$ 
as Hopf subalgebras. 
The cross-commutation relations between elements from the two tensor factors can 
be written in the form
\be{}
\label{dcr}
\langle \eta_{(1)},x_{(1)}\rangle \:\eta_{(2)}\:x_{(2)}=
\langle \eta_{(2)},x_{(2)}\rangle \:x_{(1)}\:\eta_{(1)},
\ee{}
for any $x\in \Ha$, $\eta\in \Ha^{*op}$.
In the finite dimensional case, they are equivalent to 
the Yang-Baxter equation on the canonical element $\sum_{i}e_i\tp e^i$,
where  $\{e_i\}\subset \Ha$ and $\{e^i\}\subset \Ha^{*op}$ are dual bases.
Then $\D(\Ha)$ is dual to the twisted tensor product 
\twistF{$\Ha^*$}{$\Ha_{op}$}, where 
$\F=\sum_i e_i\tp e^i$ is considered as an element from  $\Ha_{op}\tp \Ha^*$.
\begin{propn}
\label{maps}
Let $\Ha$ be a quasitriangular Hopf algebra.
The coproduct
\be{}
\Ha &\stackrel{\Delta}{\longrightarrow}&\mbox{\twist{$\Ha$}{$\Ha$}}
\label{map}
\ee{}
and the composite map
\be{}
\Ha^{*op} \stackrel{\Delta}{\longrightarrow}&\Ha^{*op}\tp\Ha^{*op}&
\stackrel{\Ru^+\tp\Ru^-}{\longrightarrow} \mbox{\twist{$\Ha$}{$\Ha$}}
\label{*map}
\ee{}
are Hopf algebra homomorphisms.
\end{propn}
\begin{proof}
Straightforward.
\end{proof}
Homomorphisms (\ref{map}) and (\ref{*map}) may be extended to $\D(\Ha)$.
For finite dimensional Hopf algebras, a proof of this statement can be found in \cite{Mj}. 
Infinite dimensional Hopf algebras like quantum groups are of the primary
interest for this article, and we present here a proof suitable for the 
general case.
\begin{thm}
\label{double}
Let $\Ha$ be a quasitriangular Hopf algebra.
Then, maps (\ref{map}) and  (\ref{*map}) define a Hopf
homomorphism 
$
\D(\Ha) \to \mbox{\twist{$\Ha$}{$\Ha$}}.
$
\end{thm}
\begin{proof}
As a linear space, the double coincides with the tensor
product of $\Ha$ and $\Ha^{*op}$, which are
embedded in $\D(\Ha)$ as Hopf subalgebras. Taking Proposition \ref{maps} into
account, it suffices to show that the permutation relations between elements
of $\Ha$ and $\Ha^{*op}$ are respected.
Applying maps (\ref{map}) and (\ref{*map}) to the both sides of identity 
(\ref{dcr}), 
we come to the equation
$$
\langle\eta\;,\Ru^-_1 \Ru^+_1 x_{(1)}\rangle \; \Ru^+_2 x_{(2)} \tp \Ru^-_2 x_{(3)}= 
\langle\eta\;,x_{(3)} \Ru^-_1 \Ru^+_1\rangle \; x_{(1)}\Ru^+_2 \tp x_{(2)}\Ru^-_2 
$$
that must hold for any $x\in \Ha$ and $\eta\in \Ha^{*op}$.
It is fulfilled indeed, because both of the elements $\Ru^\pm$ 
satisfy equation (\ref{flip}).
\end{proof}
\section{Modules over Hopf algebras.}
\label{sMHA}
\subsection{Module algebras.}
This subsection contains some facts about the modules over a Hopf algebra $\Ha$. 
An associative algebra $\A$ is called a left $\Ha$-{\em module algebra} if the multiplication
$\A\tp\A\to \A$ is a homomorphism of $\Ha$-modules. Similarly, one can consider
right modules over $\Ha$. Explicitly, for any $x\in \Ha$ and  $a,b\in \A$ the consistency
conditions read
\be{}
x\tr(a b) = \bigl(x_{(1)}\tr a\bigr) \bigl(x_{(2)}\tr b\bigr),
& &
(a b)\tl x= \bigl(a\tl x_{(1)}\bigr) \bigl(b\tl x_{(2)}\bigr)
\label{lrations}
\\
1\tr a = a, \quad x\tr 1_A = \ve(x)1_A,
& &
a\tl 1= a, \quad 1_A\tl x= \ve(x)1_A.
\label{unitlr}
\end{eqnarray}
for the left and right actions $\tr$ and $\tl$.
If $\A$ is simultaneously a left and right
module and the two actions commute, 
\be{}
\label{l-r}
x_1\tr (a \tl x_2)= (x_1\tr a) \tl x_2 ,
\quad x_1,x_2\in \Ha,\quad a\in \A,
\ee{}
then it is called {\em bimodule}. $\A$ is an $\Ha$-{\em bimodule
algebra} if its bimodule and algebra structures are 
consistent in the sense of  (\ref{lrations}--\ref{unitlr}).
\begin{example}[Adjoint action] 
A Hopf algebra $\Ha$ is a left and right module algebra over itself with respect
to the left and right adjoint actions 
\be{}
\ad(x)\tr y = x_{(1)} y \gm(x_{(2)}), \quad
y\tl\ad(x) = \gm(x_{(1)}) y x_{(2)}, 
\ee{}
for $x, y\in \Ha.$
\end{example}
\begin{example}[Dual Hopf algebra $\Ha^*$] 
A Hopf algebra $\Ha$ is a bimodule over itself with respect
to the regular actions by multiplication from the left and 
right. However, these actions do not respect the 
multiplication  in $\Ha$. On the contrary, the 
dual (coregular) actions are consistent with the multiplication
in $\Ha^*$, so the latter is an $\Ha$-bimodule algebra.
Explicitly, the coregular actions can be expressed via the 
coproduct in $\Ha^*$ and the 
pairing $\langle \:.\:,\:.\:\rangle$ between $\Ha^*$ and $\Ha$ :
\be{}
\label{lrc}
x\tr a = a_{(1)}\langle a_{(2)},x\rangle, \quad
a\tl x  = \langle a_{(1)},x\rangle a_{(2)}, 
\ee{}
where $x\in \Ha$ and $a\in \Ha^*$. 
\end{example}
\begin{example}
\label{bimod}
Let $\Ha$ be a Hopf algebra and $\A$ its bimodule algebra.
Then $\A$ is  a left  $\Ha_{op}\tp \Ha$-module algebra
with the action 
\be{}
(x\tp y)\tr a = y\tr a\tl x, \quad x, y \in\Ha_{op}\tp \Ha,\>a\in\A.
\ee{}
It is is a left  $\Ha^{op}\tp \Ha$-module algebra
with the action 
\be{}
\label{ttpaction}
(x\tp y)\tr a = y\tr a\tl \gm(x), \quad x, y \in\Ha^{op}\tp \Ha,\>a\in\A.
\ee{}
This example means that we may consider
only left modules, instead of bimodules.
\end{example}

Twist of Hopf algebras induces a transformation of 
their module algebras. 
Let $\A$ be an $\Ha$-module algebra with the multiplication $m$ and
let $\tilde \Ha\stackrel{\F}{\sim}\Ha$. The new associative multiplication 
\be{}
\label{twmod}
\tilde m (a\tp b) = m (\F_1\tr a\tp \F_2\tr b), \quad a, b \in \A,
\ee{}
can be introduced on $\A$. We denote this algebra 
by $\tilde \A$. Since $\tilde\Ha\simeq \Ha$ as
associative algebras, the action of $\Ha$ on $\A$ can be 
viewed as that of $\tilde \Ha$ on $\tilde \A$. This action
is consistent with multiplication (\ref{twmod}) in $\tilde \A$ and 
twisted coproduct (\ref{tcopr}) in $\tilde \Ha$. We say that
$\tilde \A$ and $\A$ are {\em twist-equivalent} and write 
$\tilde\A\stackrel{\F}{\sim} \A$, by the analogy with the Hopf algebras.
Also, we shall omit $\F$ and write simply $\tilde \A\sim \A$
if the exact form of $\F$ is clear from the context. 

\subsection{Comodule algebras.}
Let $\Ha$ be a Hopf algebra and $\Ha^*$ its dual.
A right $\Ha^*$-comodule algebra is an associative algebra $\A$ 
endowed with a homomorphism $\delta\colon \A\to \A\tp \Ha^*$ obeying
the coassociativity constraint
\be{}
(\id\tp \Delta)\circ\delta &=& (\delta\tp \id)\circ\delta
\label{coass}
\ee{}
and the conditions
\be{}
\delta(1_\A) =1_\A\tp 1 & & (\id\tp \ve)\circ\delta         =  \id,
\label{coun}
\ee{}
where the identity map on the right-hand side assumes the isomorphism
$\A\tp \K\simeq \A$.
As for the coproduct $\Delta$, we use symbolic notation 
$\delta(x)=x_{[1]}\tp x_{(2)}$, marking the tensor component 
belonging to $\A$ with the square brackets; the subscript
of the $\Ha^*$-component is concluded in parentheses.
Every  right $\Ha^*$-comodule $\A$ is a left $\Ha$-module,
the action being defined through the pairing $\langle \:.\:,\:.\:\rangle$
between $\Ha$ and $\Ha^*$:
\be{}
\label{coaction}
x\tr a &=& a_{[1]}\langle a_{(2)}, x  \rangle, 
\quad x\in \Ha,\> a\in  \A.
\ee{}
Suppose there is a map $\delta$ from an $\Ha$-module algebra $\A$ to $\A\tp\Ha^*$ 
such that for any
$x \in  \Ha$, $a \in  \A$, and a linear functional 
$\al\in \A^*$
\be{}
\label{delta}
\al(x\tr a) = \al(a_{[1]})\langle a_{(2)}, x\rangle,
\quad\mbox{where}\quad \delta(a)=a_{[1]}\tp a_{(2)}.
\ee{}
Then  $\A$ is an $\Ha^*$-comodule algebra with the coaction $\delta$.
Note that if $\Ha$ or $\A$ are finite dimensional, equation (\ref{delta}) 
may serve as a definition of $\delta$, so the notions of 
$\Ha$-module and $\Ha^*$-comodule are equivalent. In general, the property
of being an $\Ha^*$-comodule is stronger than of being an $\Ha$-module.

Similarly to right $\Ha^*$-comodule algebras, one can 
consider left ones. They are also right $\Ha$-module algebras.

\section{FRT- and RE-type algebras.}
\label{FRT-RE-type}
The purpose of this section is to establish a twist-equivalence
between certain classes of algebras relative to a quasitriangular
Hopf algebra $\Ha$. 
Let $V$ be a right $\Ha$-module of a finite rank over $\K$
and $V^*$ its dual. We identify the space of endomorphisms
$\End(V)$ with $V^*\tp V$ and assume $V$ to be a right
$\End(V)$-module. The representation $\rho$ of $\Ha$ on
$V$ is  a homomorphism 
$\Ha\to\End(V)$. 
The image $(\rho\tp\rho)(\Ru)\in\End^{\tp2}(V)$
of the universal  R-matrix
is denoted $R$.
Given a basis $\{e_i\}$ in $V$, the elements
$e^i_j \in \End(V)$ stand for the matrix units acting on $e_i$ by
the rule $ e_l e^i_j = e_j\delta^l_i $, where $\delta^l_j$ are the
Kronecker symbols. The multiplication in $\End(V)$ is expressed according 
to $e^i_j e^l_k = \delta^l_j e^i_k$.

The space $\End(V)$ is a bimodule for $\Ha$ or a right module for
$\Ha_{op}\tp\Ha$:
\be{}
A\tl(x\tp y) = \rho(x)A\rho(y),\quad A\in \End(V),
\quad x\tp y \in \Ha_{op}\tp\Ha.
\ee{}
By duality, the space $\End^*(V)$ is endowed with the structure a
left $\Ha_{op}\tp\Ha$-module as well.
\subsection{FRT algebra.}
\label{ssFRT}
Let $\{T^i_k\}\subset \End^*(V)$ be the basis that is dual to $\{e^k_i\}$.
The associative algebra $\Ta_R$ is generated by the matrix
coefficients $\{T^i_k\} \subset \End^*(V)$ subject to the FRT relations 
\be{}
RT_1 T_2  =  T_2T_1 R,
\label{frt}
\ee{}
where $T$ is the matrix $T=\sum_{i,j} T^i_j e^j_i$. The matrix elements 
$T^i_j$ may be thought of as linear functions on  $\Ha$; they define an algebra 
homomorphism 
\be{}
\label{TtoH}
\Ta_R\to\Ha^*. 
\ee{}
\begin{propn}
Let $\rho$ be a finite dimensional representation of $\Ha$ 
and $\Ta_R$ the FRT algebra associated with $\rho$.
Then $\Ta_R$
is a $\Ha$-bimodule algebra, with the left and right actions
\be{}
x\tr T = T\rho(x),\quad   T\tl x  = \rho(x) T, \quad  x \in \Ha
\label{lra}
\ee{}
extended from $\End^*(V)$.
It is a bialgebra, with the coproduct and counit being defined as
\be{}
\Delta(T^i_j)=\sum^n_{l=1}T^l_j\tp T_l^i,
\quad
\ve(T^i_j) =\delta^i_j.
\ee{}
Composition of the coproduct with the algebra homomorphism (\ref{TtoH}) 
applied to the  (left) right tensor factor makes
$\Ta_R$ a (left) right $\Ha^*$-comodule algebra.
\end{propn}
\begin{proof}
Actions (\ref{lra}) are extended to the actions on the tensor
algebra $\T\bigl(\End^*(V)\bigr)$ leaving invariant the ideal 
generated by (\ref{frt}). Regarding the bialgebra properties
of $\Ta_R$,
the reader is referred for the proof to \cite{FRT}. The  comodule
 structure  is inherited
from the bialgebra one, so it is obviously coassociative.
It is also an algebra homomorphism, being a composition of
two homomorphisms.
\end{proof}
\noindent
Remark that the FRT relations (\ref{frt}) arose within the quantum inverse 
scattering method and was used for systematic definition of the 
quantum group duals in \cite{FRT}.

\subsection{RE algebra.}
\label{ssRE}
Another algebra of interest, $\La_R$, is defined as the quotient 
of $\T\bigl(\End^*(V)\bigr)$ by the RE relations\footnote{Although 
$\Ta_R$ and $\La_R$ are both generated by
$\End^*(V)$, it is customary to use different letters, $T$ and $L$,
to denote their matrices of generators.}
\be{}
R_{21}L_1 R L_2  =  L_2 R_{21} L_1 R,
\label{re}
\ee{}
where $L$ is the matrix $L=\sum_{i,j} L^i_j e^j_i$ whose entries form the
set of generators.
In terms of the operator $S=PR$, where $P=\sum_{i,j}e^i_j\tp e^j_i$ is the 
permutation on $V\tp V$, relations (\ref{re}) can be written as
\be{}
SL_2 S L_2  =  L_2 S L_2 S.
\label{reS}
\ee{}
\begin{propn}
Let $\rho$ be a finite dimensional representation of $\Ha$ 
and $T^k_l\in \Ha^*$ its  matrix coefficients.
Let $L^i_j$ be the generators of the algebra $\La_R$ associated with $\rho$.
Then, $\La_R$ is a left $\Ha$-module algebra with the action 
extended from the coadjoint representation in $\End^*(V)$:
\be{}
x\tr L = \rho\bigl(\gamma(x)\bigr)L\rho(x),
\quad x \in \Ha.
\label{ada}
\ee{}
It is a right $\Ha^*$-comodule algebra with respect to the coaction
\be{}
\delta(L^i_j) = \sum_{l,k} L^l_k\tp \gamma(T^k_j) T^i_l,
\label{ada-re}
\ee{}
\end{propn}
\begin{proof}
Action (\ref{ada}) is naturally extended to  $\T\bigl(\End^*(V)\bigr)$ and
preserves
relations (\ref{reS}). The coassociativity of (\ref{ada-re}) is obvious.
To prove that $\delta$ is an algebra homomorphism, one needs to employ
commutation relations $(\ref{frt})$ and $(\ref{re})$. For details,
the reader is referred to \cite{KS}.
\end{proof}
\noindent
A spectral dependent version of the RE appeared first in
\cite{Cher}. In the form of (\ref{re}), it may be found
in articles \cite{Skl,AFS} devoted to integrable models. The algebra 
$\La_R$ was studied in \cite{KSkl,KS}. Its relation to the braid 
group of a solid handlebody was 
pointed out in \cite{K}.

\subsection{FRT- and RE-type algebras: twist-equivalence.}
\label{ssRET}
\begin{definition}
Let  $\Ha$ be a quasitriangular Hopf algebra with the universal
R-matrix $\Ru$ and $\A$ its left 
module algebra.
$\A$ is called {\em quasi-commutative} if for
any $a,b\in \A$
\be{}
\label{crel}
(\Ru_2\tr b)(\Ru_1\tr a) = a b.
\ee{}
\end{definition}
\begin{definition}
We call a quasi-commutative $\Ha^{op}\tp\Ha$-module algebra an algebra of {\em FRT-type}.
Similarly, a quasi-commutative \twist{$\Ha$}{$\Ha$}-module algebra is called
an algebra of  {\em RE-type}.
\end{definition}

\begin{example}
Let $\A$ an $\Ha$-bimodule algebra for a quasitriangular
Hopf algebra $\Ha$ with the universal R-matrix $\Ru$. 
As in Example \ref{bimod}, we can think
of it as a left $\Ha^{op}\tp\Ha$-module algebra.
Let us take $\Ru^{-1}_{13}\Ru_{24}$ for the universal R-matrix 
of $\Ha^{op}\tp\Ha$.
The algebra $\A$ is of  FRT-type if and
only if for any $a,b\in \A$
\be{}
(a\tl{\Ru_1}) (b\tl{\Ru_2})  =   ({\Ru_2}\tr b) ({\Ru_1}\tr a).
\label{rtt}
\ee{}
\end{example}
\begin{example}
The Hopf dual $\Ha^*$ viewed as an $\Ha$-bimodule with
respect to the coregular actions is an FRT-type algebra.
Relation (\ref{rtt}) is a consequence of (\ref{flip}).
\end{example}
\begin{example}
\label{frt-com}
The FRT algebra $\Ta_R$ associated with a finite dimensional
representation of $\Ha$ is of FRT-type.
Indeed, it is enough to check commutation relations (\ref{crel})
on generators. Reduced to the matrix elements  $T^i_j$ generating $\Ta_R$,
condition (\ref{rtt}) turns into (\ref{frt}).
\end{example}
\begin{propn}
\label{tw-com}
Let $\Ha$ be a quasitriangular Hopf algebra and 
$\tilde \Ha\stackrel{\F}{\sim}\Ha$.
If $\A$ is a quasi-commutative $\Ha$-module algebra, 
then the twisted $\tilde \Ha$-module algebra $\tilde\A\stackrel{\F}{\sim}\A$ 
is also quasi-commutative.
\end{propn}
\begin{proof}
Condition (\ref{crel}) holds for twisted multiplication 
(\ref{twmod}) in $\tilde\A$ and R-matrix (\ref{Rtwisted}) of $\tilde\Ha$.
\end{proof}
\begin{thm}
\label{main}
An $\Ha^{op}\tp\Ha$-module algebra is twist-equivalent
to an \twist{$\Ha$}{$\Ha$}-module algebra. 
\end{thm}
\begin{proof}
Let $\Ru^{-1}_{13}\Ru_{24}$ be the universal R-matrix for  $\Ha^{op}\tp\Ha$.
The twist from  $\Ha^{op}\tp\Ha$ to  $\Ha\tp\Ha$ with the twisting cocycle
$\Ru_{13}\in (\Ha^{op}\tp\Ha)\tp (\Ha^{op}\tp\Ha)$ turns it into 
the R-matrix $\Ru^{-1}_{31}\Ru_{24}$. Further twist with the cocycle
$\Ru_{23}\in (\Ha\tp\Ha)\tp (\Ha\tp\Ha)$ transforms  $\Ha\tp\Ha$ into 
\twist{$\Ha$}{$\Ha$} with the R-matrix $\Ru_{41}^{-1}\Ru_{31}^{-1}\Ru_{24}\Ru_{23}$,
see (\ref{ttpRm}).
\end{proof}
\begin{corollary}
Let $\A$ be an $\Ha$-bimodule algebra and  $\tilde\A$ its twist-equivalent 
left \twist{$\Ha$}{$\Ha$}-module algebra. The algebra $\tilde\A$ is 
of RE-type  if and only if $\A$ is of FRT-type. Then, 
for any $a,b\in \tilde\A$
\be{}
({\Ru_{1'}}\tr a\tl{\Ru_2}) (b\tl{\Ru_1}\tl{\Ru_{2'}}) & = & 
({\Ru_{1'}}\tr {\Ru_2}\tr b ) ({\Ru_1}\tr a\tl{\Ru_{2'}}),
\label{re1}
\ee{}
where the primes distinguish different copies of $\Ru$.
\end{corollary}
\begin{proof}
By Proposition \ref{tw-com}, $\tilde\A$ is of RE-type
{\em iff} $\A$ is of FRT-type; then
$\A$ satisfy (\ref{rtt}). Twist from $\Ha^{op}\tp\Ha$ to \twist{$\Ha$}{$\Ha$} 
converts (\ref{rtt}) to (\ref{re1}).
The action of \twist{$\Ha$}{$\Ha$} on
$\tilde\A$ is expressed through the 
right and left actions of $\Ha$ on $\A$ as in (\ref{ttpaction}).
\end{proof}
\begin{propn}
Let $\rho$ be a finite dimensional representation of $\Ha$. 
The RE algebra $\La_R$ associated with 
$\rho$ is an RE-type algebra, provided
the universal R-matrix is taken as in (\ref{ttpRm}).
It is twist-equivalent to $\Ta_R$.
\end{propn}
\begin{proof}
First, let us prove that $\La_R$ is an \twist{$\Ha$}{$\Ha$}-module
algebra.
Define the left \twist{$\Ha$}{$\Ha$}-action 
on $\End^*(V)$ by the formula 
\be{}
\label{ttp-action}
(x\tp y)\tr L = \rho\bigl(\gm(x)\bigr)L\rho(y), \quad x\tp y\in \mbox{\twist{$\Ha$}{$\Ha$}},
\ee{}
where $\{L^i_j\}\subset\End^*(V)$ is the dual basis to $\{e^j_i\}$.
This action  is uniquely extended to $\T\bigl(\End^*(V)\bigr)$
to make it 
an \twist{$\Ha$}{$\Ha$}-module algebra. It is easy to see that action
(\ref{ttp-action}) respects relations (\ref{re}) thus reducing to an action on 
$\La_R$. Now we show that $\La_R$ is of RE-type.
Abstract commutation relations (\ref{crel}) are specialized to
(\ref{re1}) for quasi-commutative \twist{$\Ha$}{$\Ha$}-module algebras.  
Evaluated on the generators $L^i_j$, they turn into the reflection equation
relations (\ref{re}). 
On the other hand, as  mentioned in Example \ref{frt-com}, the FRT relations
(\ref{frt}) are the reduction of (\ref{rtt}) to the generators of $\Ta_R$.
The twist that relates identities (\ref{rtt}) and (\ref{re1}) transforms
relations (\ref{frt}) to  (\ref{re}).
\end{proof}
\subsection{RE dual $\tilde \Ha^*$ and its properties.}
\label{ssRED}
Theorem \ref{main} allows us to define the RE analog of
the algebra $\Ha^*$ for any quasitriangular Hopf
algebra. 
\begin{definition}
{\em Reflection equation dual} $\tilde \Ha^*$ to a quasitriangular Hopf
algebra $\Ha$ is an RE-type algebra, the twist of $\Ha^*$ viewed as 
the coregular left $\Ha^{op}\tp \Ha$-module
with action (\ref{ttpaction}).
\end{definition}
Let $m$ be the multiplication in  $\Ha^*$. 
Formula  (\ref{twmod}) gives the multiplication in  $\tilde \Ha^*$
which is expressed through the universal R-matrix and the 
coregular actions of $\Ha$ on $\Ha^*$:
\be{}
\label{reprod}
\tilde m(a\tp b) =
 m\bigl(\Ru_1 \tr a \tl\Ru_{1'} \tp  b \tl \gm(\Ru_2) \tl\Ru_{2'} \bigr).
\ee{}
It follows that $\tilde \Ha^*$ is isomorphic as an associative
algebra to the braided Hopf algebra introduced by Majid 
using other arguments, \cite{Mj}.
\begin{propn}
$\tilde \Ha^*$ is a left module algebra over
the Hopf algebras $\Ha$,  $\Ha^{*op}$, and the 
double $\D(\Ha)$. 
\end{propn}
\begin{proof}
An immediate corollary of Proposition \ref{maps} and Theorem \ref{double}.
\end{proof}
Let 
$(\mbox{\twist{$\Ha$}{$\Ha$}})^*$
be the dual Hopf algebra to 
\twist{$\Ha$}{$\Ha$}. As a coalgebra, it coincides with the tensor product 
$\Ha^*\tp \Ha^*$.
The multiplication in $(\mbox{\twist{$\Ha$}{$\Ha$}})^*$ is expressed 
through coregular actions (\ref{lrc}) of the universal R-matrix of
$\Ha$:
$$
(\al\tp\bt)(\eta\tp\xi)=
\al (\Ru_2\tr \eta\tl\Ru^{-1}_{2'})\tp (\Ru_1\tr \bt\tl\Ru^{-1}_{1'})\xi,
$$
for $\al\tp\bt$ and $\eta\tp\xi\in (\mbox{\twist{$\Ha$}{$\Ha$}})^*$.
\begin{propn}
$\tilde \Ha^*$ is a right $(\mbox{\twist{$\Ha$}{$\Ha$}})^*$-comodule algebra
with the coaction expressed through the coproduct and antipode in
$\Ha^*$:
\be{}
\label{H*H*comod}
\delta(\eta)&=&\eta_{(2)}\tp \bigl(\gm(\eta_{(1)})\tp \eta_{(3)}\big),
\quad \eta\in \tilde \Ha^*.
\ee{}
\end{propn}
\begin{proof}
Map  (\ref{H*H*comod}) satisfies condition (\ref{delta})
with left action (\ref{ttpaction}). 
\end{proof}
\begin{corollary}
$\tilde \Ha^*$ is a right $\Ha^*$-comodule algebra
with the coaction 
\be{}
\label{H*comod}
\delta(\eta)&=&\eta_{(2)}\tp \gm(\eta_{(1)})\eta_{(3)},
\quad \eta\in \tilde \Ha^*.
\ee{}
$\tilde \Ha^*$ is a right $\Ha_{op}$-comodule algebra,
with the coaction 
\be{}
\label{Hcomod}
\delta(\eta)&=&\eta_{(2)}\tp 
\langle\eta_{(3)},\Ru^{-1}_{1'}\rangle\langle \gm(\eta_{(1)}),\Ru_{2}\rangle
 \Ru^{-1}_{2'}\Ru_1,
\quad \eta\in \tilde \Ha^*,
\ee{}
where $\langle\:.\:,\:.\:\rangle$ is the Hopf pairing between
$\Ha^*$ and $\Ha$.
\end{corollary}
\begin{proof}
The Hopf embedding (\ref{map})
gives rise to the reversed arrow $(\mbox{\twist{$\Ha$}{$\Ha$}})^*\to \Ha^*$
acting by multiplying the tensor factors. It remains to 
apply this homomorphism to the term in (\ref{H*H*comod}) that is 
confined in parentheses, to obtain (\ref{H*comod}). Similarly, 
evaluating the map
$(\mbox{\twist{$\Ha$}{$\Ha$}})^*\to \Ha_{op}$, which is dual to the Hopf algebra homomorphism
(\ref{*map}), we come to (\ref{Hcomod}).
\end{proof}
\begin{lemma}
\label{ve}
The counit $\ve$ of $\Ha^*$ is a character
of the RE dual  $\tilde\Ha^*$ 
\end{lemma}
\begin{proof}
Evaluating $\ve$ on the twisted product (\ref{reprod})
we find, for $a,b\in \tilde \Ha^*$,
$$
\ve\circ\tilde m(a\tp b) =
\ve \circ m\bigl(\Ru_1 \tr a \tl\Ru_{1'} \tp  b \tl \gm(\Ru_2) \tl\Ru_{2'} \bigr)=
\langle a,\Ru_{1'}\Ru_1\rangle  \langle b,\gm(\Ru_2) \Ru_{2'}\rangle =
\ve(a)\ve(b).
$$
Here we used the identity $\Ru_{1'}\Ru_1 \tp \gm(\Ru_2) \Ru_{2'}=1\tp 1$
following from (\ref{urm-antipode}).
\end{proof}
\begin{propn}
Let $\Ha$ be a quasitriangular Hopf algebra with the
universal R-matrix $\Ru$.
Consider the element $\Q=\Ru_{21}\Ru\in \Ha^{\tp 2}$.
The map $\tilde \Ha^* \to \Ha$ defined by the correspondence
\be{}
\label{Qmap}
\eta\to \langle\eta,\Q_1\rangle \Q_2, \quad \eta\in \tilde \Ha^*,
\ee{}
where $\langle\:.\:,\:.\:\rangle$ is the Hopf pairing between
$\Ha^*$ and $\Ha$,
is a homomorphism of the coadjoint and adjoint $\Ha$-module
algebras.
\end{propn}
\begin{proof}
According to Lemma \ref{ve}, the counit of $\Ha^*$ is
a character of $\tilde \Ha^*$. 
Applying the antipode to the right tensor factor of
(\ref{Hcomod}) we pass to a left comodule structure
with respect to the Hopf algebra $\Ha^{op}$. 
Applying the counit to the left tensor factor belonging
to $\tilde \Ha^*$, we obtain an equivariant homomorphism
of $\Ha$-module algebras $\tilde \Ha^*\to \Ha$.
Taking composition of $(\ve\tp\gm)$ with coaction 
(\ref{Hcomod}) we obtain (\ref{Qmap}):
\be{}
(\ve\tp \gm)\circ\delta(\eta)&=&
\langle\eta_{(2)},\Ru^{-1}_{1'}\rangle\langle \gm(\eta_{(1)}),\Ru_{2}\rangle
\gm(\Ru^{-1}_{2'}\Ru_1)
\n\\&=&
\langle\eta,\gm(\Ru_{2})\Ru^{-1}_{1'}\rangle
\gm(\Ru_1) \gm(\Ru^{-1}_{2'})
=
\langle\eta,\Ru_{2}\Ru_{1'}\rangle
\Ru_1\Ru_{2'}.
\ee{}
Here we used $(\gm\tp\gm)(\Ru)=\Ru$ and $(\id\tp\gm)(\Ru^{-1})=\Ru$,
see (\ref{urm-antipode}).
Equivariance of map (\ref{Qmap}) may be derived from the comodule structure
(\ref{Hcomod}). However, it is readily seen directly.
Since $\Delta(x)\Q=\Q\Delta(x)$ for every $x\in \Ha$, one has 
$\gm(x_{(1)}) \Q_1 x_{(2)}\tp \Q_2 = \Q_1\tp  x_{(1)}\Q_2 \gm(x_{(2)})$.
Then,
\be{}
\langle x_{(2)}\tr \eta \tl \gm(x_{(1)}), \Q_1\rangle \;\Q_2 =
\langle  \eta , \gm(x_{(1)}) \Q_1 x_{(2)}\rangle \;\Q_2 =
\langle  \eta ,  \Q_1\rangle  \; x_{(1)}\Q_2 \gm(x_{(2)}),
\ee{}
for $x\in \Ha$ and $\eta\in \tilde \Ha^*$.
\end{proof}
\section{Applications: equivariant quantization on $G$-spaces.}
\label{sAQGS}
In this section, we specialize the constructions of the previous sections to the deformation 
quantization situation when $\Ha=\U_h(\g)$, a quantum group corresponding
to a semisimple Lie algebra $\g$.
\subsection{Quantization on the group space $G$.}
\label{ssQGS}
Let $G$ be a semisimple Lie group equipped and $\g$ its Lie algebra.
An element $\xi\in \g$ generates left and right 
invariant vector fields
$$
(\xi^l \tr f)(g)=\frac{d}{dt}f(ge^{t\xi})|_{t=0},
\quad 
(f\tl\xi^r)(g) =\frac{d}{dt}f(e^{t\xi}g)|_{t=0},
$$
defining the left and right actions of the algebra $\U(\g)$ on
functions on $G$. Given an element $\psi\in \U(\g)$, by $\psi^r$ and $\psi^l$
we  correspondingly denote its extensions by the right- and left-invariant
differential operators on  $G$. We use notation 
$\xi^{\ad}=\xi^l-\xi^r$ for the vector field generated by
the element $\xi\in \g$ via the adjoint action 
$a\to g^{-1} a g$ of  $G$ on itself.

Let $r\in \wedge^2\g$ be a classical r-matrix and $\omega$ the invariant
symmetric element such that $r+\omega$ satisfies the classical Yang-Baxter-Equation,
\cite{Dr1}. Let  $\U_h(\g)$ be the corresponding quantum group. 
\begin{propn}
The RE dual $\tilde \U^*_h(\g)$ to the quantum group 
$\U_h(\g)$ is a \twist{$\U_h(\g)$}{$\U_h(\g)$}-equivariant
quantization of the Poisson bracket 
\be{}
r^{r,r}+r^{l,l}-r^{l,r}-r^{r,l} + (\omega^{r,l}-\omega^{l,r})
\label{rebr}
\ee{}
on the group $G$.
It turns into a $\U_h(\g)$-equivariant quantization of the bracket
\be{}
r^{\ad,\ad}+ (\omega^{r,l}-\omega^{l,r})
\label{rebr1}
\ee{}
via the Hopf algebra embedding $\Delta\colon\U_h(\g)\to \mbox{\twist{$\U_h(\g)$}{$\U_h(\g)$}}$.
\end{propn}
\begin{proof}
The element 
$r^-=(- r) \oplus r$  from the exterior square $\wedge^{2}(\g\oplus \g)$
generates the bivector field $r^-_G=-r^{r,r}+r^{l,l}$ on $G$ via the action
of $\U(\g)_{op}\tp \U(\g)$; $r^-_G$ coincides with the Drinfeld-Sklyanin bracket
on $G$.
The algebra $\U^*_h(\g)$ is the $\U_h(\g)^{op}\tp \U_h(\g)$-equivariant quantization 
of $r^-_G$, \cite{Ta}.
The twist with the cocycle $\Ru_{13}=1+h(r_{13}+ \omega_{13}) +o(h)$, converts $\U^*_h(\g)$ 
into the $\U_h(\g)\tp \U_h(\g)$-equivariant quantization of the bracket 
$r^{r,r}+r^{l,l}$. Indeed, at the infinitesimal level this procedure adds 
the term $2r^{r,r}$ to the bracket $r^-_G$.
Further twist with the cocycle
$\Ru_{23}=1+h(r_{23}+ \omega_{23}) +o(h)$ 
leads to the algebra $\tilde \U^*_h(\g)$. In terms of Poisson brackets, this 
operation adds the term $-(r^{l,r}+\omega^{l,r})-(r^{r,l}-\omega^{r,l})$ 
to $r^{r,r}+r^{l,l}$ thus resulting in (\ref{rebr}).
The last statement of the proposition is straightforward.
\end{proof}
\subsection{Quantization of polynomial functions on matrices.}
\label{ssQM}
Let  $V$ be a complex vector space and $\rho$ a homomorphism of  $\U(\g)$ into the matrix 
algebra $\M=\End(V)$. It induces a bimodule structure
on $\M^*$:
$$
(x\tr f)(A)= f\bigl(A \rho (x)\bigr), \quad  (f\tl x)(A)  = f\bigl(\rho (x) A\bigr),
$$
for $x\in \U(\g)$, $f\in \M^*$, and $A\in \M$.
Let $\Omega\in \M^{\tp2}$ be the image of the invariant symmetric element
$\omega$. 
Denote by $\M^{\Omega}$ the cone of matrices 
$$\M^{\Omega}=\{ A\in \M|\; [\Omega , A\tp A]=0\}.$$ 
Evidently, $\M^{\Omega}$ 
is an algebraic variety, it is closed under the matrix multiplication and invariant 
with respect to the two-sided action of the group $G$.

\begin{remark}
We would like to stress that we do not restrict the consideration to fundamental 
representations of $\U(\g)$. The subspace $\M^\Omega$ coincides with $\M=\End(V)$ 
only for $\g=sl(n,\C)$ and  $V=\C^n$.
\end{remark}
\begin{propn}
The quotient of the algebra $\La_R$ by the torsion is a 
\twist{$\U_h(\g)$}{$\U_h(\g)$}-equivariant quantization of the Poisson bracket
\be{}
r^{r,r}+r^{l,l}-r^{l,r}-r^{r,l} + (\omega^{r,l}-\omega^{l,r})
\label{rebr3}
\ee{}
\end{propn}
\begin{proof}
It is easy to see that $\M^{\Omega}$ is the maximal subspace in 
$\M$ where bracket (\ref{rebr3}) is Poisson. It was proven in
\cite{DS} that the algebra $\Ta_R$ is a 
$\U_h(\g)_{op} \tp \U_h(\g)$-equivariant quantization on  $\M^{\Omega}$.
Applying the RE twist, we obtain the algebra $\La_R$ as the quantization 
on $\M^{\Omega}$.
This twist transforms the bracket $-r^{l,l}+r^{r,r}$ on $\M^{\Omega}$ to 
bracket (\ref{rebr3}). This proves the statement.
\end{proof}
Note that bracket (\ref{rebr3}) goes over into $r^{\ad,\ad}+(\omega^{r,l}-\omega^{l,r})$
after restriction of $\g\oplus \g$ to the diagonal subalgebra.
Then  $\La_R$ becomes an equivariant quantization of this bracket
with respect to $\U_h(\g)\subset\mbox{\twist{$\U_h(\g)$}{$\U_h(\g)$}}$.

\bigskip
e-mail: donin@macs.biu.ac.il\\
e-mail: mudrova@macs.biu.ac.il; 
\end{document}